\documentclass[a4paper,12pt]{amsart}

\usepackage{float}
\usepackage{graphicx}

\usepackage{amsmath}
\usepackage[mathscr]{eucal}

\textwidth=\paperwidth \advance\textwidth by-2in
\def\A{\mathcal{A}}
\def\V{\mathcal{V}}
\def\Z{\mathbf{Z}}

\def\R{\mathbf{R}}

\def\lk{\operatorname{lk}}
\def\Ker{\operatorname{Ker}}
\def\Coker{\operatorname{Coker}}

\def\GL{\operatorname{GL}}

\def\spmatrix#1{\left[
  \begin{smallmatrix} #1 \end{smallmatrix}
  \right]\ignorespaces}
\theoremstyle{plain}
\newtheorem{thm}{Theorem}[section]

\newtheorem{cor}[thm]{Corollary}
\newtheorem{lem}[thm]{Lemma}
\theoremstyle{definition}
\newtheorem*{exmp}{Example}
\newtheorem{remark}[thm]{Remark}

\begin{document}

\title[Fibred knots and twisted alexander invariants]
{Fibred knots and twisted alexander invariants}
\date{June 26, 2002 (First Edition: July 18, 2001)}
\author{Jae Choon Cha}
\address{Department of Mathematics,
Indiana University,
Bloomington, IN 47405, U.S.A}
\email{jccha\char`\@indiana.edu}
\def\subjclassname{\textup{2000} Mathematics Subject Classification}
\expandafter\let\csname subjclassname@1991\endcsname=\subjclassname
\expandafter\let\csname subjclassname@2000\endcsname=\subjclassname
\subjclass{57M25}
\keywords{Fibred Knots, Twisted Alexander Invariants}

\begin{abstract}
We introduce a new algebraic topological technique to detect
non-fibred knots in the three sphere using the twisted Alexander
invariants.  As an application, we show that for any Seifert matrix of
a knot with a nontrivial Alexander polynomial, there exist infinitely
many non-fibered knots with the given Seifert matrix.  We illustrate
examples of knots that have trivial Alexander polynomials but do not
have twisted Alexander invariants of fibred knots.
\end{abstract}

\maketitle

\section{Introduction and main results}
\label{sec:intr-main-results}

A knot in the three sphere is called a \emph{fibred knot} if its
exterior is the total space of a fibre bundle over the circle.  It is
well known that the Alexander polynomial of a fibred knot is monic,
that is, the coefficients of the highest and the lowest degree terms
are units ($\pm 1$ in~$\Z$).  In general, the converse is not true,
although it is true for several kinds of knots, including 2-bridge
knots, closures of positive braids, and knots up to ten crossings.
Some further necessary conditions for fibred knots have been known;
Gabai developed a geometric procedure to detect
fibredness~\cite{Gabai:1986-1}.  Silver and Williams found an entropy
invariant which vanishes for fibred
knots~\cite{Silver-Williams:1999-1}.

In this paper we introduce a new algebraic topological technique to
detect non-fibred knots, using the \emph{twisted Alexander
invariants}.  Since a twisted version of the Alexander polynomial of a
knot was first defined by Lin~\cite{Lin:1990-1}, the twisted Alexander
invariants have been studied by several authors, including
Wada~\cite{Wada:1994-1}, Jiang and Wang~\cite{Jiang-Wang:1993-1}, and
Kitano~\cite{Kitano:1996-1}.  Especially, in a remarkable work of Kirk
and Livingston~\cite{Kirk-Livingston:1999-2, Kirk-Livingston:1999-3},
a topological definition of the twisted Alexander polynomial was
introduced.  For a surjection of the fundamental group $\pi$ of a
complex $X$ onto $\Z$ and a representation of $\pi$ over a field~$F$,
they defined the twisted Alexander module to be a specific twisted
homology group of~$X$, which is a module over the principal ideal
domain~$F[\Z]$, and defined the twisted Alexander polynomial to be a
polynomial representing the order of (the torsion part of) the twisted
Alexader module.  Some known results on the Alexander polynomial was
extended to the twisted case.  In particular, a slicing obstruction
was obtained from the twisted Alexander polynomials associated to
specific representations of the fundamental groups of prime power fold
cyclic branched covers, by relating them with the Casson-Gordon
invariants.

We generalize the approach of Kirk and Livingston for a noetherian
unique factorization domain~$R$ which is \emph{not} necessarily a
field.  Roughly speaking, first we define the twisted Alexander ideal
to be the elementary ideal of the twisted Alexander module, and then
we define the twisted Alexander polynomial to be the greatest common
divisor of specific generators of the elementary ideal.  When $R$ is a
field, the twisted Alexander polynomial of ours coincides with that of
Kirk and Livingston if the former is nontrivial.  (For details, see
Section~\ref{sec:twisted-alex-poly}.)

Using the above terminologies, our main result is stated as follows.
Suppose that every submodule of a free module of finite rank over the
base domain $R$ is again free of finite rank.  The main example to
keep in mind is $R=\Z$.

\begin{thm}
\label{thm:nonfibred-knot-twisted-alex-inv}
If $K$ is a fibred knot and $\rho$ is a representation of the
fundamental group of a cyclic cover of $S^3$ branched along $K$ that
factors through a finite group, then the twisted Alexander invariants
associated to $\rho$ have the following properties:
\begin{enumerate}
\item
The twisted Alexander module is $R[\Z]$-torsion.
\item
The twisted Alexander ideal is a principal ideal generated by the
twisted Alexander polynomial.
\item
The twisted Alexander polynomial is monic.
\end{enumerate}
\end{thm}

This can be viewed as a natural generalization of the property of the
classical Alexander polynomial of a fibred knot.  We remark that if
$R$ is a field, then the conclusions of
Theorem~\ref{thm:nonfibred-knot-twisted-alex-inv} hold for any knot
with a nontrivial twisted Alexander polynomial.  This is the reason
why we work with the twisted Alexander invariants over a domain $R$
which is not a field.

In Section~\ref{sec:fibred-knot}, we develop a practical method to
compute the twisted Alexander invariants of a fibred knot from (the
given representation and) the homotopy type of a \emph{monodromy} on a
fibre surface, which can be described as a map of a graph, or
equivalently as an automorphism of a free group.
Theorem~\ref{thm:nonfibred-knot-twisted-alex-inv} is proved using this
computational method.

As an application of
Theorem~\ref{thm:nonfibred-knot-twisted-alex-inv}, we prove a
realization theorem of Seifert matrices by non-fibred knots.

\begin{thm}
\label{thm:nonfibred-knot-with-given-seifert-surface}
If $A$ is a Seifert matrix of a knot with a nontrivial Alexander
polynomial, there exist infinitely many non-fibred knots with Seifert
matrix~$A$.
\end{thm}

In Section~\ref{sec:non-fibred-knot-with-given-seifert-matrix}, we
describe explicitly how to construct non-fibred knots for any given
Seifert matrix with a nontrivial Alexander polynomial.  Non-fibredness
is proved using the twisted Alexander invariants associated to
``abelian'' representations (whose images are abelian subgroups
of~$\GL_n$).

Since the Alexander polynomial of a nontrivial fibred knot is
nontrivial, it follows that there are infinitely many non-fibred knots
which are indistinguishable from a given nontrivial fibred knot via
invariants derived from Seifert matrices, including Alexander modules,
torsion invariants, signatures, and Blanchfield linking forms.  In a
recent work of Cochran that appeared subsequent to ours, this result
has been generalized for higher order Alexander
modules~\cite{Cochran:2002-1}.

In Section~\ref{trivial-alex-poly-example}, we illustrate by examples
that there are knots that have \emph{trivial} Alexander polynomials
but do not have twisted Alexander invariants of fibred knots, i.e.,\ do
not satisfy the conclusions of
Theorem~\ref{thm:nonfibred-knot-twisted-alex-inv}.  For this purpose
we use non-abelian representations.

\section{Twisted Alexander invariants}
\label{sec:twisted-alex-poly}

We begin with the definitions of the twisted Alexander invariants,
which generalize the twisted Alexander polynomial defined
in~\cite{Kirk-Livingston:1999-2}.  Throughout this section, we assume
that $R$ is a noetherian unique factorization domain.  Let $X$ be a
finite CW-complex and $\tilde X$ be its universal covering space.
$\pi_1(X)$~acts on the left of $\tilde X$ as the covering
transformation group.  Let $\epsilon\colon \pi_1(X)\to \langle
s\rangle$ be a surjection, where $\langle s\rangle$ denotes the
infinite cyclic group generated by~$s$.  We identify the group ring
$R[\langle s\rangle]$ with the Laurent polynomial ring $R[s,s^{-1}]$
and denote it by~$\Lambda$.  Suppose that $\pi_1(X)$ acts on the right
of a free $R$-module~$V$ of finite rank via a representation $\pi_1(X)
\to \GL(V)$.  Then $\Lambda \otimes_R V$ becomes a
$\Lambda$-$R[\pi_1(X)]$ bimodule under the actions given by
$s^k\cdot(s^n\otimes v)=s^{n+k}\otimes v$ and $(s^n\otimes v)\cdot g =
s^n \epsilon(g) \otimes v\cdot g$ for $v \in V$ and $g \in \pi_1(X)$.
Let $C_*(\tilde X;R)$ be the cellular chain complex of $\tilde X$ with
coefficient~$R$, which is a left $R[\pi_1(X)]$-module.  The twisted
cellular complex of $X$ with coefficient $\Lambda\otimes V$ is defined
to be the following chain complex of left $\Lambda$-modules:
$$
C_*(X;\Lambda\otimes V) = (\Lambda\otimes V) \otimes_{R[\pi_1(X)]}
C_*(\tilde X;R).
$$
The twisted homology $H_*(X;\Lambda \otimes V)$ is defined to be the
homology of~$C_*(X;\Lambda \otimes V)$.  We call $H_1(X; \Lambda
\otimes V)$ the \emph{twisted Alexander module}.  For notational
convenience, we denote it by~$A$.

Since $V$ is finitely generated, so is $A$ as a $\Lambda$-module.
Since $R$ is noetherian, so is~$\Lambda$, and hence $A$ is a finitely
presentable $\Lambda$-module.  Choose a presentation of~$A$, and let
$n$ and $m$ be the numbers of generators and relations, respectively.
Let $P$ be the $n\times m$ matrix associated to the presentation,
i.e.,\ the $(i,j)$-entry of $P$ is the coefficient of the $i$-th
generator in the $j$-th relation.  Let $D$ be the set of the
determinants of all $n\times n$ submatrices obtained by removing
$(m-n)$ columns from~$P$.  The ideal $\A$ in $\Lambda$ generated by
$D$ is called the \emph{elementary ideal} of~$A$.  It is known that
the elementary ideal is an invariant of the $\Lambda$-module~$A$,
which is independent of the choice of~$P$.  For a proof,
see~\cite[p.101]{Crowell-Fox:1977-1}.  A similar argument shows that
the greatest common divisor $\Delta$ of elements of $D$ is also an
invariant of~$A$, which is well defined up to multiplication of $us^n$
with $u$ a unit in~$R$.  (If $n>m$, $A=\{0\}$ and $\Delta=0$ by a
convention.)  We call $\A$ and $\Delta$ the \emph{twisted Alexander
ideal} and the \emph{twisted Alexander polynomial}, respectively.

\begin{remark}
\label{rmk:elem-ideal-annihilator}
The elementary ideal of a module is contained in the annihilator
ideal.  In particular, elements of $\A$ annihilate~$A$.
\end{remark}

\begin{remark}
When $R$ is a field and $A$ is $\Lambda$-torsion, our definition is
equivalent to that of Kirk and
Livingston~\cite[$\mathsection$2]{Kirk-Livingston:1999-2}.  In this
case $\Lambda$ is a principal ideal domain, and by the classification
theorem of finitely generated modules over~$\Lambda$, $A$ is
decomposed into a direct sum of cyclic modules $\bigoplus_i
R[s,s^{-1}]/\langle d_i \rangle$, where $d_i$ is an element
of~$\Lambda$.  Kirk and Livingston defined the twisted Alexander
polynomial to be the product of all nonzero~$d_i$.  It is equal to the
order of the torsion part of~$A$.  Since the diagonal matrix with
diagonals $d_i$ is a presentation matrix of~$A$, $\Delta$~is equal to
the product of all~$d_i$.  This shows that the two definitions
coincide if $\Delta\ne 0$, or equivalently $A$ is torsion.
Furthermore, when $R$ is a field, $\A$~is the principal ideal~$\langle
\Delta \rangle$.  In general, $\A$~is contained
in~$\langle\Delta\rangle$, however, the converse is not true if $R$ is
not necessarily a field.
\end{remark}

\begin{remark}
Using Fox's calculus, one can compute the boundary map $C_2(X;\Lambda
\otimes V) \to C_1(X;\Lambda \otimes V)$ from a presentation
of~$\pi_1(X)$.  Consequently the twisted Alexander invariants can also
be computed.  In~\cite{Kirk-Livingston:1999-2}, this method was used
as the main computational technique and was used to relate this
topological version of twisted Alexander polynomial with Wada's
invariants~\cite{Wada:1994-1}.  Because this method is not used in
this paper, we do not proceed into further details.
\end{remark}

Specifically, we define the twisted Alexander invariants of an
oriented knot $K$ in $S^3$ as follows.  The first homology of the
exterior $E$, which is obtained by removing an open tubular
neighborhood of $K$ from~$S^3$, is an infinite cyclic group generated
by an element $t$ such that $\lk(K,t)=+1$.  Let $N$ be the $d$-fold
cyclic cover of~$E$.  The image of the composition $\pi_1(N) \to
\pi_1(E) \to H_1(E)=\langle t \rangle$ is the subgroup generated by
$s=t^d$.  Thus it induces a surjection $\epsilon \colon \pi_1(N) \to
\langle s \rangle$ so that the twisted Alexander invariants of $N$ are
defined for any representation of~$\pi_1(N)$.  In this paper, we will
consider only representations that factor through the fundamental
group of the $d$-fold cyclic cover $M$ of $S^3$ branched along $K$.
View $N$ as a subspace of~$M$, and let $i_*\colon \pi_1(N)\to \pi_1(M)$
be the homomorphism induced by the inclusion.  For a representation
$\rho$ of $\pi_1(M)$, we denote the twisted Alexander module, ideal
and polynomial of $N$ associated to $\epsilon$ and $\rho i_*$ by
$A_K^\rho$, $\A_K^\rho$ and $\Delta_K^\rho$, respectively.

\begin{remark}
The twisted Alexander invariants of the exterior $E$ associated to
$\pi_1(E)\to H_1(E)=\langle t \rangle$ and a representation
of~$\pi_1(E)$ are also useful in studying knots.  This version appears
in some literature including~\cite{Wada:1994-1, Kitano:1996-1,
Kirk-Livingston:1999-2, Kirk-Livingston:1999-3}, where the last two
concern our version (with field coefficients) as well.
\end{remark}

\section{Fibred knots}
\label{sec:fibred-knot}

Let $K$ be a fibred knot, and let $M$ be the $d$-fold branched cover
of $K$ as before.  In this section, we are interested in a special
case of representations that factor through finite groups; we assume
that a given representation $\rho$ of $\pi_1(M)$ is decomposed~as
$$
\rho\colon \pi_1(M) \xrightarrow{\phi} G \to \GL(V)
$$
where $\phi$ is a homomorphism into a finite group $G$ and $G\to
\GL(V)$ is a representation in a free $R$-module~$V$ of finite rank.
Without any loss of generality, we may assume that $\phi$ is
surjective.  In addition, we assume that $R$ has the property that
every submodule of a free $R$-module of finite rank is again free of
finite rank.  For example, the ring $\Z$ has this property.

We will compute the twisted Alexander invariants of $K$ associated to
$\rho$ from a monodromy of~$K$.  Let $F$ be a fibre surface of $K$ and
$h\colon F\to F$ be a monodromy such that
$$
E= \R\times F /(r,x)\sim (r+1, h(x)), \quad r\in \R,\ x\in F
$$
is the exterior of $K$ and $\{0\}\times \partial F$ represents a
preferred longitude of~$K$.  Then the $d$-fold cyclic cover $N$ of the
exterior is given by
$$
N = \R\times F/(r,x) \sim (r+d,h^d(x)), \quad r\in \R,\ x\in F.
$$
The preferred generator of the covering transformation group acts on
$N$ by $[r,x] \mapsto [r+1,h(x)]$.  We note that any fibre surface $F$
is connected, and so is~$N$.

First of all, we need to compute $\pi_1(N)$ and~$\pi_1(M)$.  Since $N$
can be viewed as a quotient space of $[0,d]\times F$ under an obvious
identification, $\pi_1(N)$ is expressed as an HNN-extension of
$\pi_1(F)$.  An explicit description is as follows.  Fix a basepoint
$*$ on~$\partial F$.  Let $\gamma\colon [0,1]\to \partial F$ be a path
from $*$ to $h(*)$ such that the loop $\tau(t) = [t,\gamma(t)]$ in
$E=\R \times F/\mathord{\sim}$ is a positive meridian of~$K$.
Identifying $\pi_1(F,h(*))$ with $\pi_1(F)=\pi_1(F,*)$ under the
isomorphism induced by~$\gamma$, $h$ induces an endomorphism $h_*$
of~$\pi_1(F)$ which is given by $h_*([\delta])=[\gamma\cdot h\delta
\cdot \gamma^{-1}]$.  Then $\pi_1(N)$ is presented as
$$
\pi_1(N) = \langle s, \pi_1(F) \mid szs^{-1}=h_*^d(z) \text{ for }
z\in \pi_1(F) \rangle,
$$
and the map $\epsilon$ defined in the previous section is equal to the
surjection $\pi_1(N)\to \pi_1(N)/\pi_1(F)=\langle s \rangle$.  Since
$\pi_1(M)$ is obtained from $\pi_1(N)$ by killing~$s$, $\pi_1(M)$ is
isomorphic to the quotient group of $\pi_1(F)$ modulo the normal
subgroup generated by $\{z^{-1}h_*^d(z)\mid z\in \pi_1(F)\}$.

Let $\tilde F$ be the connected regular covering of $F$ associated to
the composition
$$
\alpha\colon \pi_1(F) \to \pi_1(M) \xrightarrow{\phi} G,
$$
that is, the kernel of $\alpha$ is equal to the image of the injection
$\pi_1(\tilde F) \to \pi_1(F)$ induced by the covering projection.
$G$~acts on $\tilde F$ as the covering transformation group.  We need
the following lemma, which is an easy exercise in the covering space
theory.  Since the author has not found a proof in the literature, we
give a proof for completeness.

\begin{lem}
Suppose that $p\colon (\tilde X,\tilde x_0) \to (X,x_0)$ is a regular
covering projection, $\tilde X$ and $X$ are connected and locally path
connected, and $f \colon X \to X$ be a map.  Choose a path $\gamma$
from $x_0$ to~$f(x_0)$ and let $f_*$ be the endomorphism on $\pi_1(X)$
defined by $f_*([\delta])=[\gamma \cdot f\delta \cdot
\gamma^{-1}]$. If $z^{-1} f_*(z) \in p_*\pi_1(\tilde X)$ for all $z\in
\pi_1(X)$, then $f$ is lifted to a map $\tilde f\colon \tilde X\to
\tilde X$ which commutes with the action of covering transformations,
i.e., $\tilde f\tau=\tau\tilde f$ for any covering
transformation~$\tau$.
\end{lem}
\begin{proof}
Choose $\tilde x_0 \in p^{-1}(x_0)$.  Let $\tilde \gamma$ be the lift
of $\gamma$ starting at~$\tilde x_0$, and let $\tilde x_1$ be the
endpoint of~$\tilde \gamma$.  By the hypothesis $f_*(z)\in z\cdot
p_*\pi_1(\tilde X)$, $f_*p_*\pi_1(\tilde X) \subset p_*\pi_1(\tilde
X)$.  Thus there is a lift $\tilde f$ of $f$ such that $f(\tilde
x_0)=\tilde x_1$ by the lifting criterion.  Let $\tau$ be a covering
transformation on~$\tilde X$.  Since $\tilde f \tau$ and $\tau \tilde
f$ are lifts of the same map $fp$, it suffices to show that $\tilde f
\tau(\tilde x_0)$ and $\tau \tilde f(\tilde x_0)$ coincide by the
uniqueness of a lift.  Let $\tilde\delta$ be a path from $\tilde x_0$
to $\tau(\tilde x_0)$.  Then $\tilde f \tau(\tilde x_0)$ is the
endpoint of $\tilde f \tilde \delta$.  Since $\tau$ is the
transformation associated to the loop $\delta=p\tilde \delta$,
$\tau(\tilde x_1)$~is the endpoint of the path $\tilde \gamma^{-1}
\cdot \tilde\delta \cdot \tau\tilde \gamma$.  Since $[\delta^{-1}
\cdot \gamma \cdot f\delta \cdot \gamma^{-1}] = [\delta]^{-1}
f_*([\delta]) \in p_*\pi_1(\tilde X)$, the endpoints of $\tilde f
\tilde \delta$ and $\tilde \gamma^{-1} \cdot \tilde\delta \cdot
\tau\tilde \gamma$ are the same.
\end{proof}

In our case, $z^{-1}h^d_*(z)$ is contained in the kernel of~$\alpha$
for any $z \in \pi_1(F)$.  Thus by the lemma, the homeomorphism
$h^d\colon F\to F$ is lifted to a homeomorphism $\tilde h^d\colon
\tilde F\to \tilde F$ which commutes with the action of~$G$.  Now
$\langle s \rangle \oplus G$ acts on $\R\times \tilde F$ by
$(s^n,g)\cdot (r,w) = (r+nd, (\tilde h^d)^n(g\cdot w))$.  It is easily
checked that the orbit space is $N$, and the projection $\R\times
\tilde F \to N$ is a covering projection with covering transformation
group $\langle s \rangle \oplus G$.  This shows that $\R\times \tilde
F$ is the regular covering of~$N$ associated to the homomorphism
$\epsilon \oplus \phi i_* \colon \pi_1(N) \to \langle s \rangle \oplus
G$, where $i_*\colon \pi_1(N)\to \pi_1(M)$ is the map induced by the
inclusion.

Recall that the twisted Alexander module $A_K^\rho$ is defined to be
the twisted homology group $H_1(N; \Lambda\otimes V)$ where
$\Lambda=R[\langle s \rangle]$.  It is equal to the first homology of
the chain complex
$$
(\Lambda\otimes V) \otimes_{R[\langle s \rangle \oplus G]}
C_*(\R\times \tilde F;R) \cong V\otimes_{R[G]} C_*(\R\times \tilde
F;R).
$$
As an $R$-module, it can be viewed as the twisted homology group
$H_1(\R\times F;V) = H_1(F;V)$.  Since $F$ is a connected surface with
nonempty boundary, $F$ has the homotopy type of a graph (1-complex)
with one vertex.  Therefore $H_1(F;V)$ can be computed from a chain
complex
$$
\cdots \to 0 \to V\otimes_{R[G]} R[G]^n \xrightarrow{\partial_1}
V\otimes_{R[G]} R[G] \to 0
$$
where $n$ is the number of edges of the graph.  Since $H_1(F;V) =
\Ker(\partial_1)$ is a submodule of $V\otimes_{R[G]} R[G]^n = V^n$, it
is a free $R$-module of finite rank.

The action of $s$ on $H_1(F;V)=H_1(V\otimes_{R[G]} C_*(\tilde F))$ is
given by the homomorphism induced by $\tilde h^d \colon \tilde F\to
\tilde F$.  Let $H$ be a matrix associated to the induced map by
choosing an $R$-basis of~$H_1(F;V)$.  Then $sI-H$ is a presentation
matrix of $A_K^\rho$, as a $\Lambda$-module, where $I$ is the identity
matrix.  Since it is a square matrix, $\A_K^\rho$ is the principal
ideal generated by $\Delta_K^\rho(s)=\det(sI-H)$.  The coefficient of
the highest term of $\Delta_K^\rho(s)$ is equal to $\det(I)=1$.  The
constant term of $\Delta_K^\rho(s)$ is~$\det(H)$, which is a unit in
$R$ since $\tilde h^d$ is an homeomorphism on~$\tilde F$.

We summarize the above discussion as a theorem.

\begin{thm}
\label{thm:twisted-alex-inv-of-fibred-knot}
\begin{enumerate}
\item
$A_K^\rho$ is presented by the matrix $sI-H$, as a $\Lambda$-module.
\item
$\A_K^\rho$ is the principal ideal generated by~$\Delta_K^\rho(s)$.
\item
$\Delta_K^\rho(s)=\det(sI-H)$ is a monic polynomial.
\end{enumerate}
\end{thm}

The last two conclusions of
Theorem~\ref{thm:nonfibred-knot-twisted-alex-inv} are immediate
consequences of Theorem~\ref{thm:twisted-alex-inv-of-fibred-knot}.  By
Remark~\ref{rmk:elem-ideal-annihilator}, $A_K^\rho$ is annihilated by
$\Delta_K^\rho(s)$, and so the first conclusion of
Theorem~\ref{thm:nonfibred-knot-twisted-alex-inv} follows.  The
following consequence of our discussion will be useful later.

\begin{cor}
\label{cor:annhilator-of-twisted-alex-module}
$A_K^\rho$ is annihilated by a monic polynomial.
\end{cor}

\begin{remark}
If $V=R[G]$ and $G\to \GL(V)$ is the regular representation,
Theorem~\ref{thm:twisted-alex-inv-of-fibred-knot} and
Corollary~\ref{cor:annhilator-of-twisted-alex-module} are true
\emph{without} the assumption that every submodule of a free
$R$-module of finite rank is free of finite rank.  For, in this case,
the twisted homology $H_1(F;V)$ is equal to $H_1(\tilde F;R)$, and
hence it is always a free $R$-module of finite rank since $\tilde F$
is a surface without closed components.
\end{remark}

We finish this section with an example illustrating the above
computational method for fibred knots. 

\begin{exmp}
Let $K$ be the trefoil knot.  There is a well known fibre structure of
the exterior of~$K$, e.g.,\ see~\cite{Rolfsen:1976-1} for a detailed
description.  We need only the following fact: a~monodromy of $K$ has
the homotopy type of a map $h$ on a graph $B$ with one vertex and two
oriented edges $x$ and~$y$, which is defined by $h(x)=y^{-1}$ and
$h(y)=xy$.

Then the fundamental group of the double branched cover $M$ of $K$ is
given by
$$
\pi_1(M)=\langle x,y \mid x=h^2(x)=y^{-1}x^{-1}, y=h^2(y)=y^{-1}xy
\rangle.
$$
By simplifying relations, $\pi_1(M)$ is a cyclic group of order $3$
generated by $x=y$.  Let $\rho$ be the regular representation of
$\pi_1(M)$ over~$\Z$.

We will compute the twisted Alexander invariants associated to~$\rho$.
By the above discussion, $A_K^\rho=H_1(B; \Z[\langle s\rangle]\otimes
\Z[\pi_1(M)])=H_1(\tilde B)$ where $\tilde B$ is the regular cover of
$B$ associated to the homomorphism $\pi_1(B)\to \Z_3$ given by
$x,y\mapsto 1$.  Obviously $\tilde B$ is again a graph; $\tilde B$ has
$3$ vertices $v_0,v_1,v_2$ and $6$ edges $x_0$, $x_1$, $x_2$, $y_0$,
$y_1$, $y_2$ where $\partial x_i=\partial y_i=v_{i+1}-v_i$ (indices
are modulo~$3$), and the covering projection $\tilde B\to B$ sends
$x_i$ and $y_i$ to $x$ and $y$, respectively.  $\pi_1(\tilde B)$~can
be identified with the free subgroup of $\pi_1(B)$ generated by
$a=xy^{-1}$, $b=xax^{-1}$, $c=x^2ax^{-2}$ and $d=x^3$, and hence
$H_1(\tilde B)$ is the free abelian group generated by (the homology
classes of) $a$, $b$, $c$ and $d$.  The action of $s$ on $H_1(\tilde
B)$ is easily computed by evaluating the values of $h^2$ on $a$, $b$,
$c$ and~$d$; for example,
$h^2(a)=y^{-1}x^{-1}y^{-1}x^{-1}y=d^{-1}cad^{-1}c^{-1}d$
in~$\pi_1(\tilde B)\subset \pi_1(B)$, and by abelianizing, $s\cdot
a=a-d$ in~$H_1(\tilde B)$.  By computing the action on the other
generators in a similar way, we obtain a matrix
$$
H=\begin{bmatrix}
1 & 0 & 0 & 1 \\
0 & 1 & 0 & 1 \\
0 & 0 & 1 & 1 \\
-1 & -1 & -1 & 2
\end{bmatrix}
$$
which represents the action of $s$ on~$A_K^\rho$.  Thus $A_K^\rho$ is
presented by~$sI-H$, $\Delta_K^\rho(s)=\det(sI-H)=s^4-s^3-s+1$, and
$\A_K^\rho$ is the principal ideal generated by~$\Delta_K^\rho(s)$.
\end{exmp}

\section{Non-fibred knots with given Seifert matrices}
\label{sec:non-fibred-knot-with-given-seifert-matrix}

In this section we prove
Theorem~\ref{thm:nonfibred-knot-with-given-seifert-surface}.  We begin
with a lemma, which is a consequence of well known classical results.

\begin{lem}
\label{lem:alex-poly-and-homology-of-branched-cover}
The Alexander polynomial $\Delta_{K_0}(t)$ of a knot $K_0$ is nontrivial
if and only if the first homology group of the $d$-fold cyclic cover
of $S^3$ branched along $K_0$ is nontrivial for some~$d$.
\end{lem}
\begin{proof}
The order of the first homology of the $d$-fold cyclic branched cover
is given by the resultant
$$
R_d = \Big| \prod_{s=0}^{d-1} \Delta_{K_0}(e^{2\pi is/d}) \Big|
$$
where $R_d=0$ if the homology is an infinite
group~\cite{Fox:1956-1,Gordon:1978-1,Burde-Zieschang:1985-1}. Thus if
$\Delta_{K_0}(t)$ is trivial, $R_d=1$ for all~$d$.

Conversely, if $\Delta_{K_0}(t)$ is nontrivial, the equation
$\Delta_{K_0}(t)=0$ has a nonzero complex root~$w$.  If $w$ is a
$d$-th root of unity, then $R_d=0$.  In~\cite{Gordon:1972-1}, Gordon
proved that if $w$ is not a root of unity, then the nonzero values of
$R_d$ are unbounded.  (Actually, more is known; Riley proved that the
nonzero values of $R_d$ grow exponentially
in~$d$~\cite{Riley:1990-1}.)  This completes the proof.
\end{proof}

By the lemma, it suffices to prove

\begin{thm}
\label{thm:realization-by-nonfibred-knot}
If the first homology of the $d$-fold cyclic branched cover of a knot
$K_0$ is nontrivial for some~$d$, and $A$ is a Seifert matrix of $K_0$,
then there exist infinitely many non-fibred knots with Seifert
matrix~$A$.
\end{thm}

The remaining part of this section is devoted to the proof of
Theorem~\ref{thm:realization-by-nonfibred-knot}.  Choose $d$ such that
the $d$-fold cyclic branched cover $M_0$ of a given knot $K_0$ has
nontrivial first homology.  Choose a surjection $\chi_0$ of $H_1(M_0)$
onto a cyclic group $\Z_r$ of order $r\ge 2$, and choose a Seifert
surface $F$ of~$K$.  For the given data $K_0$, $\chi_0$ and~$F$, we
will construct non-fibred knots that admit the same Seifert form as
that of $K_0$ defined on~$F$.

\begin{figure}[hb]
\begin{center}
\includegraphics{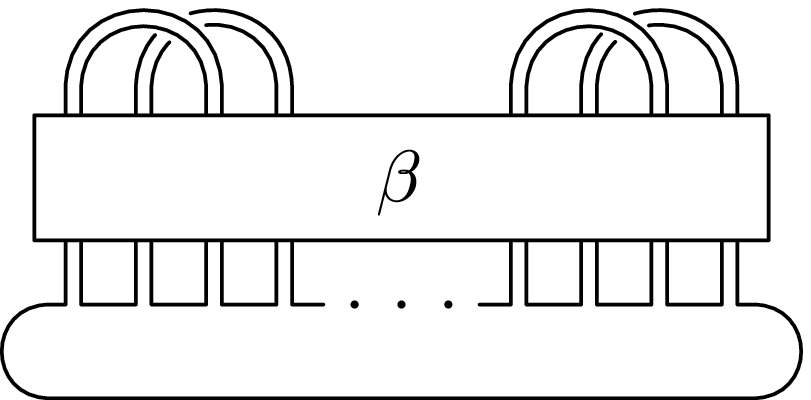}
\end{center}
\caption{}
\label{fig:seifert-surface}
\end{figure}

\begin{figure}[hb]
\begin{center}
\includegraphics{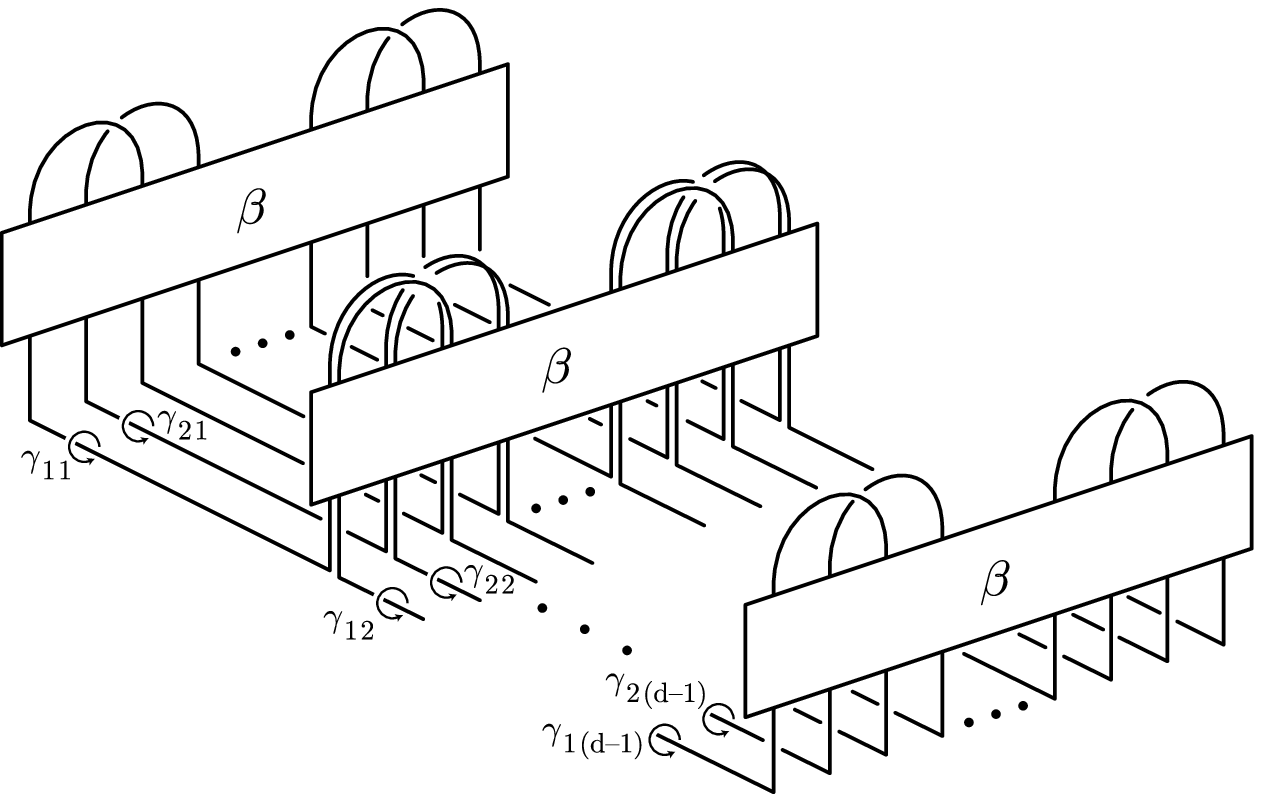}
\end{center}
\caption{}
\label{fig:branched-cover}
\end{figure}

We may assume that $F$ is a handlebody with one $0$-handle and $s$
1-handles, and by an isotopy we may assume that $F$ is embedded in
$S^3$ as in Figure~\ref{fig:seifert-surface}, where $\beta$ is a
framed $(2s)$-string link.  By a method of Akbulut and
Kirby~\cite{Akbulut-Kirby:1979-1}, $M_0$ is obtained by surgery on a
$s(d-1)$-component link $L$ shown in Figure~\ref{fig:branched-cover}.
Furthermore, a presentation of $H_1(M_0)$ is obtained as follows.
Denote meridians of components of $L$ by $\gamma_{ij}$ ($1 \le i \le
s$, $1 \le j \le d-1$) as in Figure~\ref{fig:branched-cover}, and let
$A$ be the Seifert matrix of $F$ with respect to the generators of
$H_1(F)$ represented by the 1-handles.  Then $H_1(M_0)$ is generated
by $\{\gamma_{ij}\}$ and
$$
\begin{bmatrix}
A+A^T  & -A^T              \\
-A     & A+A^T  & \ddots \\
       & \ddots & \ddots & -A^T \\
       &        &  -A    & A+A^T
\end{bmatrix}_{s(d-1)\times s(d-1)}
$$
is a presentation matrix of~$H_1(M_0)$ with respect
to~$\{\gamma_{ij}\}$.

We claim that $\chi_0(\gamma_{ij}) \ne \chi_0(\gamma_{i(j+1)})$ in
$\Z_r$ for some $i$ and~$j$.  Suppose not.  Then $\chi_0(\gamma_{ij})$
is independent of~$j$; let $x_i=\chi_0(\gamma_{ij})$.  Since all
relations on $\gamma_{ij}$ are killed by~$\chi_0$, we have
$x(A+A^T)-xA=-xA^T+x(A+A^T)=0$ where $x$ is the row vector with
entries~$x_i$.  From this we easily obtain $x(A-A^T)=0$.  Since $A$ is
a Seifert matrix of a knot, $A-A^T$ is nonsingular and $x=0$.  This
implies that $\chi_0$ is a trivial map, a contradiction.  Thus the
claim is true.

\begin{figure}[hb]
\begin{center}
\includegraphics{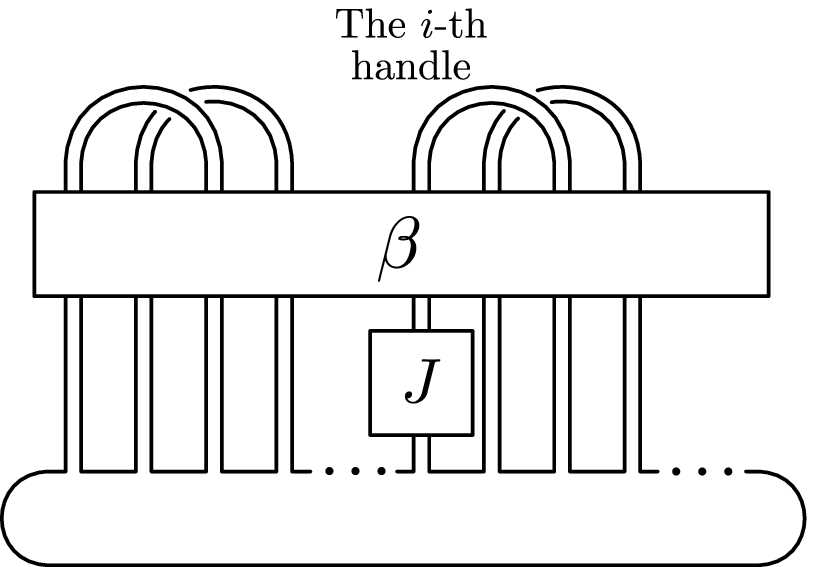}
\end{center}
\caption{}
\label{fig:satellite-knot}
\end{figure}

Let $K$ be the knot obtained by tying a knot $J$ along the $i$-th
handle of $F$ so that the Seifert form is unchanged.  See
Figure~\ref{fig:satellite-knot}.  Let $n$ be the order of
$\chi_0(\gamma_{ij}) - \chi_0(\gamma_{i(j+1)})$ in $\Z_r$.  Note that
$n\ge 2$.  Let $M_J$ be the $n$-fold cyclic cover of $S^3$ branched
along~$J$.  Using the twisted Alexander invariants, we will show

\begin{thm}
\label{thm:nonfibredness-of-satellite}
If $H_1(M_J)$ is nontrivial, then $K$ is not a fibred knot.
\end{thm}

Actually, there exist infinitely many knots $J$ satisfying the
hypothesis of Theorem~\ref{thm:nonfibredness-of-satellite}. For
example, the following lemma says that the hypothesis is true for any
knot $J$ having the same Seifert matrix as that of the figure eight
knot.

\begin{lem}
\label{lem:branched-cover-homology-realization}
if $\spmatrix{1 & 1 \\ 0 & -1}$ is a Seifert matrix of~$J$, then the
first homology of the $n$-fold cyclic cover of $S^3$ branched along
$J$ is nontrivial for all $n\ge 2$.
\end{lem}
\begin{proof}
Since $A$ is nonsingular, $H^n-I$ is a presentation matrix of the
first homology of the $n$-fold cyclic branched cover of~$J$, where
$H=A^{-1}A^T = \spmatrix{2 & -1 \\ -1 & 1}$.  It suffices to show that
$H^n-I$ is not unimodular for $n\ge 2$.  Let $H^n=\spmatrix{a_n & b_n
\\ b_n & c_n}$.  (Note that $H$ is a symmetric matrix.)  Then it is
easily shown that $a_n \ge 5$, $b_n \le -3$ and $c_n \ge 2$ for $n\ge
2$, by an induction.  Since $\det(H)=1$, we have
$$
\det(H^n-I)=(a_n-1)(c_n-1)-b_n^2 = \det(H^n)+1-a_n-c_n = 2-a_n-c_n \le -5.
$$
\end{proof}

Since a different choice of $J$ produces a different knot $K$ (e.g.\
by the uniqueness of the torus decomposition of knot complements),
Theorem~\ref{thm:realization-by-nonfibred-knot} follows
Theorem~\ref{thm:nonfibredness-of-satellite}.

\begin{proof}[Proof of Theorem~\ref{thm:nonfibredness-of-satellite}]

Let $M$ be the $d$-fold cyclic branched cover of~$K$.  Since the tying
operation does not change the Seifert matrix, the Akblute-Kirby method
gives the same presentations of $H_1(M)$ and $H_1(M_0)$, and hence
there is a natural isomorphism between them.  Let $\chi\colon
H_1(M)\cong H_1(M_0) \xrightarrow{\chi_0} \Z_r$ and let
$$
\rho\colon \pi_1(M) \to H_1(M) \xrightarrow{\chi} \Z_r \to \GL(\Z[\Z_r])
$$
where the last map is the regular representation of~$\Z_r$.  In order
to show that $K$ is not fibred, we investigate the twisted Alexander
invariants associated to $\rho$.  Let $\Lambda=\Z[\langle s \rangle]$
and $N$ be the $d$-fold cyclic cover of the exterior of $K$ as before.
By definition, $A_K^\rho$ is equal to the twisted homology $H_1(N;\V)$
where $\V=\Lambda\otimes_\Z \Z[\Z_r]=\Z[\langle s \rangle \oplus
\Z_r]$.

Let $U$ be an unknotted solid torus in $S^3-F$ which links the $i$-th
handle of $F$.  The exterior of $K$ is obtained from the exterior of
$K_0$ by removing the interior of $U$ and by filling with the exterior
of~$J$ along the boundary.  The meridian (resp.\ the longitude) of $J$
is identified with a curve on $\partial U$ which is homotopic to the
core of $U$ (resp.\ null-homotopic in~$U$).  Since the linking number
of $U$ and $K_0$ is zero, $U$~is lifted to the $d$-fold cyclic cover
$N_0$ of the exterior of~$K_0$.  $N$~is obtained by removing the
interiors of all lifts of $U$ from $N_0$ and filling $d$ copies of the
exterior of $J$ along the boundaries.  Viewing $N_0$ as a subspace
of~$M_0$, $\gamma_{ij}-\gamma_{i(j+1)}$ is homologous to the core of a
lift of $U$ in~$M_0$.  In~$N$, the boundary of that lift bounds a copy
of the exterior of~$J$.  Denote it by~$E_J$.

Let $Y$ be the closure of~$N-E_J$.  Applying the Mayer-Vieotoris
theorem to $N = E_J \cup Y$, we obtain an exact sequence
\begin{multline*}
\cdots \to H_1(\partial E_J;\V) \to H_1(E_J;\V) \oplus H_1(Y;\V) \to
H_1(N;\V) \\
\to H_0(\partial E_J;\V) \to H_0(E_J;\V)\oplus H_0(E_J;\V).
\end{multline*}

The twisted homologies of $E_J$ and $\partial E_J$ have a simple
structure as follows.  First we observe that (1) the linking number of
$U$ and $K_0$ is zero, and (2) the map $\chi$ sends the meridian of
$J$ in $E_J$ to the element $\chi_0(\gamma_{ij}-\gamma_{i(j+1)})$,
which is of order $n$ in~$\Z_r$.  From the observations, the $(\langle
s \rangle \oplus \Z_r)$-covering $\widetilde{E}_J$ of $E_J$ is a union
of infinitely many copies of the $n$-fold cyclic cover $N_J$ of~$E_J$
and we have $H_*(E_J;\V)=H_*(\V \otimes_{\Z[\langle s \rangle \oplus
\Z_r]}C_*(\widetilde{E}_J)) = H_*(\widetilde
E_J)=H_*(N_J)^{r/n}\otimes_{\Z} \Lambda$.  Similarly, $H_*(\partial
E_J;\V)=H_*(\partial N_J)^{r/n} \otimes_{\Z} \Lambda$.  Therefore
$H_0(\partial E_J;\V) \to H_0(E_J;\V)$ is an isomorphism, and
$\Coker\{H_1(\partial E_J;\V) \to H_1(E_J;\V)\}$ is isomorphic to
$\Coker\{H_1(\partial N_J) \to H_1(N_J)\}^{r/n} \otimes_\Z \Lambda =
H_1(M_J)^{r/n}\otimes_\Z \Lambda$.  By the below lemma, an annihilator
of $H_1(N;\V)$ annihilates $H_1(M_J)\otimes_\Z \Lambda$ as well.

\begin{lem}
Suppose that $A$ and $B$ are modules and $C$ is a submodule of
$A\oplus B$.  Let $p\colon A\oplus B \to A$ be the canonical
projection.  Then an annihilator of $(A\oplus B)/C$ annhilates
$A/p(C)$ as well.
\end{lem}
\begin{proof}
Viewing $A$ and $B$ as submodules of $A\oplus B$, $(A\oplus B)/C =
((A+C)+(B+C))/C = (A+C)/C+(B+C)/C$.  Thus $(A+C)/A \cong A/A\cap C$ is
a submodule of $(A\oplus B)/C$.  Since $p(C)$ contains $A\cap C$,
$A/p(C)$ is a quotient of $A/(A\cap C)$.  The conclusion follows.
\end{proof}

Since $H_1(M_J)$ is nontrivial, $H_1(M_J)\otimes\Lambda$ is never
annihilated by any nonzero monic polynomial, and hence so
is~$A_K^\rho$.  By
Corollary~\ref{cor:annhilator-of-twisted-alex-module}, $K$ is not a
fibred knot.
\end{proof}

\begin{remark}
Our construction is similar to one in~\cite{Kirk-Livingston:1999-2},
which was used to illustrate that the twisted Alexander module of the
\emph{complement} of a knot is not necessarily $\Lambda$-torsion.  The
above argument shows that an analogous result holds for the twisted
Alexander module of a \emph{branched cover} of a knot as well; if
$H_1(M_J)$ is free abelian, $A_K^\rho$ is not $\Lambda$-torsion, since
any annihilator of $A_K^\rho$ annihilates $\Lambda$ and hence must be
zero by the proof of Theorem~\ref{thm:nonfibredness-of-satellite}.
This is a significant difference between the twisted Alexander
invariants and the classical Alexander invariants.
\end{remark}

\begin{remark}
If a knot $K$ has a trivial Alexander polynomial, twisted Alexander
invariants associated to abelian representations are none more than
the classical ones; indeed all abelian representations are trivial
since the first homologies of cyclic branched covers of $K$ always
vanish by Lemma~\ref{lem:alex-poly-and-homology-of-branched-cover}.
However, invariants associated to non-abelian representations are
still interesting, as shown in the next section.
\end{remark}

\section{Examples with trivial Alexander polynomials}
\label{trivial-alex-poly-example}

In this section we illustrate examples of knots that have trivial
Alexander polynomials but do not have twisted Alexander modules of
fibred knots.

\begin{figure}[hb]
\begin{center}
\includegraphics{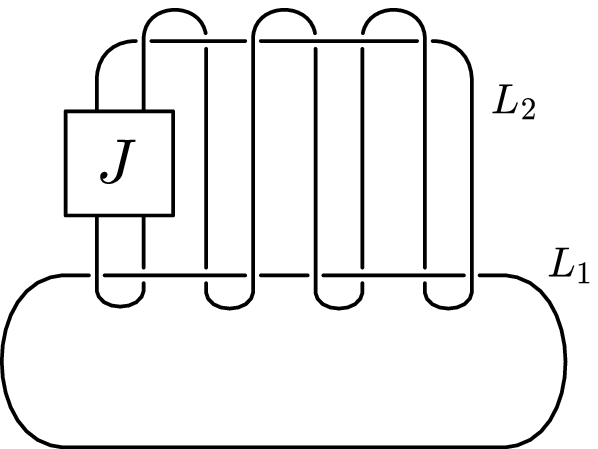}
\end{center}
\caption{}
\label{fig:triv-alex-poly}
\end{figure}

For a knot~$J$, consider the link shown in
Figure~\ref{fig:triv-alex-poly}.  Denote its components by $L_1$ and
$L_2$ as in the figure.  As usual, $J$ is tied so that the writhe of
the diagram is unchanged.  Performing $(1/1)$-surgery along~$L_2$, the
ambient space still remains $S^3$ but the other component $L_1$
becomes a knotted circle~$K$.

We remark that a similar construction was used to produce a knot with
a given Alexander polynomial in ~\cite{Levine:1965-2}
and~\cite[7.C.5]{Rolfsen:1976-1}.  Actually, by the same arguments, it
is easily seen that for any $J$ our construction produces a knot $K$
with a trivial Alexander polynomial.

Let $M_J$ be the 3-fold cyclic cover of $S^3$ branched along~$J$.  Our
goal is to prove

\begin{thm}
\label{thm:triv-alex-poly-nonfibred-knot}
If $H_1(M_J)$ is nontrivial, $K$ does not have twisted Alexander invariants of fibred knots.
\end{thm}

\begin{proof}
First we consider a special case where $J$ is unknotted.  Let $K_0$ be
the knot obtained from an unknot $J$ by the above construction.  Let
$M_0$ be the double branched cover of~$K_0$.  Cutting $S^3$ along the
obvious disk bounded by~$L_1$ and pasting two copies, we obtain a
surgery diagram of~$M_0$ with two components, which is shown in
Figure~\ref{fig:triv-alex-poly-branched-cover}.

\begin{figure}[ht]
\begin{center}
\includegraphics{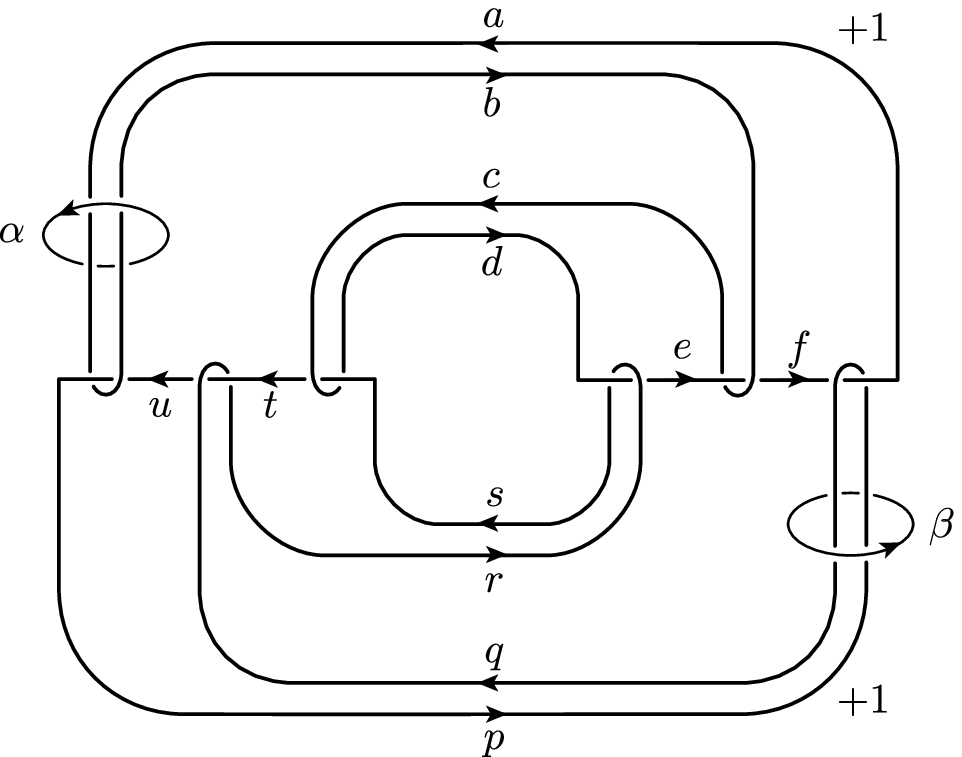}
\end{center}
\caption{}
\label{fig:triv-alex-poly-branched-cover}
\end{figure}

The Wirtinger presentation of the diagram in
Figure~\ref{fig:triv-alex-poly-branched-cover} (forgetting the
framing) is as follows: there are $12$ generators $a$, $b$, $c$, $d$,
$e$, $f$, $p$, $q$, $r$, $s$, $t$, $u$, and $12$ relations
\begin{align*}
a&=q^{-1}fq, & b&=p^{-1}ap, & c&=e^{-1}be, & d&=scs^{-1}, & e&=rdr^{-1}, &
f&=b^{-1}eb,\\
p&=b^{-1}ub, & q&=a^{-1}pa, & r&=t^{-1}qt, & s&=drd^{-1}, & t&=csc^{-1}, &
u&=q^{-1}tq
\end{align*}
where any one of the relations is redundant.  Adding the relations
$$
qpes^{-1}r^{-1}bf^{-1}=1, \quad batd^{-1}c^{-1}qu^{-1}=1
$$
which represent the effect of surgery, we obtain a presentation
of~$\pi_1(M_0)$.

We define a homomorphism $\phi_0$ of $\pi_1(M_0)$ into~$A_5$, the
group of even permutations on $\{1,\ldots,5\}$, by assigning values to
generators as follows.  (Cycle notations are used to represent
elements of~$A_5$.)
\begin{align*}
\phi_0(a)&=(132),&
\phi_0(b)&=(142),&
\phi_0(c)&=(125),\\
\phi_0(d)&=(243),&
\phi_0(e)&=(145),&
\phi_0(f)&=(152),\\
\phi_0(p)&=(13542),&
\phi_0(q)&=(15432),&
\phi_0(r)&=(12534),\\
\phi_0(s)&=(14523),&
\phi_0(t)&=(15324),&
\phi_0(u)&=(14352).
\end{align*}
It is tedious but straightforward to verify that all relations are
killed by~$\phi_0$, and in addition, $\phi_0$~is a surjection.

A representation for the general case is induced by the homomorphism
$\phi_0$ as follows.  Let $K$ be the knot obtained from a knot $J$
(not necessarily unknotted) by the above construction.  Let $U$ and
$V$ be tubular neighborhoods of the curves $\alpha$ and $\beta$ shown
in Figure~\ref{fig:triv-alex-poly-branched-cover}.  Then the double
branched cover $M$ of $K$ is obtained from $M_0$ by removing the
interiors of $U$ and $V$ and by filling two copies $E_J$ and $E_J'$ of
the exterior of $J$ along $\partial U$ and $\partial V$, respectively.
By the Seifert-van Kampen theorem, $\pi_1(M)$ is an amalgamated
product of $\pi_1(E_J)$, $\pi_1(E_J')$ and~$\pi_1(M-(E_J\cup E_J'))$.
Let $\phi_1$ and $\phi_2$ be homomorphisms of $H_1(E_J)$ and
$H_1(E_J')$ into $A_5$ which send the meridian to $(243)$ and $(253)$,
respectively.  Since $\phi_0(\alpha)=\phi_0(ab^{-1})=(243)$ and
$\phi_0(\beta)=\phi_0(pq^{-1})=(253)$, the homomorphisms
\begin{gather*}
\pi_1(M-(E_J\cup E_J')) \to \pi_1(M_0) \xrightarrow{\phi_0} A_5 \\
\pi_1(E_J) \to H_1(E_J) \xrightarrow{\phi_1} A_5 \\
\pi_1(E_J') \to H_1(E_J') \xrightarrow{\phi_2} A_5
\end{gather*}
induce a homomorphism $\phi\colon \pi_1(M) \to A_5$.  Let $\rho$ be
the representation of $\pi_1(M)$ obtained by composing $\phi$ with the
regular representation of~$A_5$.

Now we are ready to apply the arguments of the previous section.
Since $\phi$ sends the meridian of $J$ in $E_J$ to an element $(243)$,
that is of order 3 in $A_5$, an annihilator of $A_K^\rho$ annihilates
$H_1(M_J)\otimes_\Z \Lambda$ as in the proof of
Theorem~\ref{thm:nonfibredness-of-satellite}.  If $H_1(M_J)$ is
nontrivial, $H_1(M_J)\otimes_\Z \Lambda$ is never annihilated by any
monic polynomial, and so is~$A_K^\rho$.  Therefore $K$ does not have
twisted Alexander invariants of fibred knots by
Corollary~\ref{cor:annhilator-of-twisted-alex-module}.
\end{proof}

\begin{remark}
By Lemma~\ref{lem:branched-cover-homology-realization}, there are
infinitely many knots satisfying the hypothesis of
Theorem~\ref{thm:triv-alex-poly-nonfibred-knot}.
\end{remark}
\begin{remark}
Even in the case where $H_1(M_J)$ is trivial (in particular when $J$
is unknotted), $K$ is not fibred.  Indeed, since $\pi_1(M)$ is
nontrivial, $K$ is a nontrivial knot having a trivial Alexander
polynomial.
\end{remark}

\bibliographystyle{amsplainabbrv}
\bibliography{research}

\providecommand{\bysame}{\leavevmode\hbox to3em{\hrulefill}\thinspace}
\begin{thebibliography}{10}

\bibitem{Akbulut-Kirby:1979-1}
S.~Akbulut and R.~Kirby, \emph{Branched covers of surfaces in $4$-manifolds},
  Math. Ann. \textbf{252} (1979/80), no.~2, 111--131.

\bibitem{Burde-Zieschang:1985-1}
G.~Burde and H.~Zieschang, \emph{Knots}, Walter de Gruyter \& Co., Berlin,
  1985.

\bibitem{Cochran:2002-1}
T.~D. Cochran, \emph{{Noncommutative Knot Theory}}, arXiv:math.GT/0206258.

\bibitem{Crowell-Fox:1977-1}
R.~H. Crowell and R.~H. Fox, \emph{Introduction to knot theory},
  Springer-Verlag, New York, 1977, Reprint of the 1963 original, Graduate Texts
  in Mathematics, No. 57.

\bibitem{Fox:1956-1}
R.~H. Fox, \emph{Free differential calculus. {I}{I}{I}. {S}ubgroups}, Ann. of
  Math. (2) \textbf{64} (1956), 407--419.

\bibitem{Gabai:1986-1}
D.~Gabai, \emph{Detecting fibred links in ${S}\sp 3$}, Comment. Math. Helv.
  \textbf{61} (1986), no.~4, 519--555.

\bibitem{Gordon:1972-1}
C.~Gordon, \emph{Knots whose branched cyclic coverings have periodic homology},
  Trans. Amer. Math. Soc. \textbf{168} (1972), 357--370.

\bibitem{Gordon:1978-1}
\bysame, \emph{Some aspects of classical knot theory}, Knot theory (Proc. Sem.,
  Plans-sur-Bex, 1977), Springer, Berlin, 1978, pp.~1--60.

\bibitem{Jiang-Wang:1993-1}
B.~J. Jiang and S.~C. Wang, \emph{Twisted topological invariants associated
  with representations}, Topics in knot theory (Erzurum, 1992), Kluwer Acad.
  Publ., Dordrecht, 1993, pp.~211--227.

\bibitem{Kirk-Livingston:1999-2}
P.~Kirk and C.~Livingston, \emph{Twisted {A}lexander invariants, {R}eidemeister
  torsion, and {C}asson-{G}ordon invariants}, Topology \textbf{38} (1999),
  no.~3, 635--661.

\bibitem{Kirk-Livingston:1999-3}
\bysame, \emph{Twisted knot polynomials: inversion, mutation and concordance},
  Topology \textbf{38} (1999), no.~3, 663--671.

\bibitem{Kitano:1996-1}
T.~Kitano, \emph{Twisted {A}lexander polynomial and {R}eidemeister torsion},
  Pacific J. Math. \textbf{174} (1996), no.~2, 431--442.

\bibitem{Levine:1965-2}
J.~P. Levine, \emph{A characterization of knot polynomials}, Topology
  \textbf{4} (1965), 135--141.

\bibitem{Lin:1990-1}
X.-S. Lin, \emph{Representations of knot groups and twisted alexander
  polynomials}, to appear in Acta Math.\ Sinica (Series B).

\bibitem{Riley:1990-1}
R.~Riley, \emph{Growth of order of homology of cyclic branched covers of
  knots}, Bull. London Math. Soc. \textbf{22} (1990), no.~3, 287--297.

\bibitem{Rolfsen:1976-1}
D.~Rolfsen, \emph{Knots and links}, Publish or Perish Inc., Berkeley, Calif.,
  1976, Mathematics Lecture Series, No. 7.

\bibitem{Silver-Williams:1999-1}
D.~S. Silver and S.~G. Williams, \emph{Knot invariants from symbolic dynamical
  systems}, Trans. Amer. Math. Soc. \textbf{351} (1999), no.~8, 3243--3265.

\bibitem{Wada:1994-1}
M.~Wada, \emph{Twisted {A}lexander polynomial for finitely presentable groups},
  Topology \textbf{33} (1994), no.~2, 241--256.

\end{thebibliography}

\bigskip

\end{document}